\documentclass{article}
\usepackage{amsmath}
\font\tengoth=eufm10
\font\sevengoth=eufm7
\font\fivegoth=eufm5
\newfam\gothfam
\textfont\gothfam=\tengoth \scriptfont\gothfam=\sevengoth
\scriptscriptfont\gothfam=\fivegoth

\newtheorem{theorem}{Theorem}[section]

\newtheorem{definition}[theorem]{Definition}
\newtheorem{corollary}[theorem]{Corollary}
\newtheorem{lemma}[theorem]{Lemma}

\newtheorem{remark}[theorem]{Remark}

\def\blacksquare{\hbox to .60em{\vrule width .60em height .60em}}

  \font\bb=msbm10 
\catcode`\@=12  

\def\na{\nabla}

\def\?{\'e}
\def\{\`e}
\def\?{\`a}
\def\{\`u}
\def\{\c c}
\def\hb {\hfil \break}
\def\n {\vskip 0.2cm \noindent }
\def\scirc{\,{\raise 0.8pt\hbox{$\scriptstyle\circ$}}\,}
\def\ins{\,{\raise 0.2cm \hbox{ $\scriptstyle \circ$}}\,}

\def  \é{\'e}
\def\è{\`e}
\def\à{\`a}
\def\ù{\`u}
\def\ç{\c c$\!\!\!$}

\date{}
\title{Ordinary holomorphic webs of codimension one }
\author{Vincent Cavalier \and Daniel Lehmann}
 
\begin{document}
  
 \centerline {\bf  Ordinary holomorphic webs of codimension one } 
 \bigskip
   
 \centerline  {by Vincent CAVALIER and Daniel LEHMANN.}
 
  \begin{abstract}
 
 \bigskip

 \bigskip
 
To any    $d$-web of codimension one on a holomorphic $n$-dimensional  manifold $M $ ($d>n$), we associate an analytic subset $S$ of $M$. We  call {\bf ordinary} the  webs for  which $S$ has   a dimension at most $n-1$ or is empty. This condition is generically satisfied.

 \medskip
  
We  prove that the  rank of    a ordinary    $d$-web has 
  an  upper-bound $\pi'(n,d)$           
  which, for $n\geq 3$,   is strictly smaller than the bound $\pi(n,d)$ of Castelnuovo. 
This bound is optimal.  
 \medskip
   
Setting $c(n,h)=\begin{pmatrix} n-1+h\\h\\  \end{pmatrix}$,  let $k_0$ be the integer such that $c(n,k_0)\leq d<c(n,k_0+1)$. The number  $\pi'(n,d)$ is then equal\smallskip  

- to 0 for $d<c(n,2)$, \smallskip 
 
 - and to 
   $\sum _{h= 1}^{k_0}\bigl(d-c(n,h)\bigr)$ for $d\geq  c(n,2)$.
  
  \medskip 
Moreover, if   $d$ is precisely equal to $c(n,k_0)$,    we define    a holomorphic connection       on a holomorphic bundle $\cal E$ of rank $\pi'(n,d )$,  such that the space  of abelian relations   is isomorphic to    the space  of  holomorphic sections of $\cal E$ with vanishing      covariant derivative : the curvature of this connection, which  generalizes the  Blaschke curvature, is then an obstruction for the rank of the web to reach the  value $\pi'(n,d )$.
\medskip

 When $n=2$,   any web is ordinary, $\pi'(2,d)=\pi(2,d)$, any $d$ may be written  $c(2,k_0)$, and    we recover   the  results given locally in [P][He1].

 \end{abstract}

\maketitle
\section{Introduction}

 \indent A holomorphic  $d$-web $W$ of codimension one  on a $n$-dimensional holomorphic manifold $M$ being given   ($d>n$), 
    we denote by $M_0$ the open set in $M$ where the web   is   locally defined by $d$   holomorphic foliations ${\cal F}_i$ of codimension one, all non-singular and with   tangent spaces   to the leaves distinct   at any point.

    We shall say that the web is in {\bf     strong   general position} if  \emph {any} subset of $n$ leaves among the $d$ leaves through a point of $M_0$ are in general position. If we assume only that  \emph {there exists} a  subset of $n$ leaves among the $d$ leaves through a point of $M_0$ which are in general position (but not necessarily any $n$ of them), we shall say that the web is in {\bf    weak general position}. Most of our results below will require only this weaker condition\footnote{More generally,  if there exists   $\ell$  and not more of the leaves through a point   which are  in general position ($\ell\leq n)$, there exists locally a holomorphic foliation ${\cal G}$ of  codimension $ \ell$, such that the   web is locally the pullback   of a $d$-web of codimension 1 on a   $\ell$-dimensional manifold transversal to ${\cal G}$. Then, many of the  results below   remain valid, after replacing   $n$ by $\ell$.}.

 On an open set  $U$ in $M_0$ sufficiently small for the web to be  defined as above by the data of $d$ distinct local foliations ${\cal F}_i$ (such an open set is called ``open set of distinguishability"), an abelian relation  is the data   of a family $(\omega_i)_i$ of   holomorphic 1-forms $\omega_i,\  (1\leq i\leq d)$ such that 
 
$(i)$ each $\omega_i$ is closed,
 
$(ii)$ the vector fields tangent to    the local foliation ${\cal F}_i$  belong to  the kernel of    $\omega_i$,
 
$(iii)$ the sum $\sum_{i=1}^d \omega_i$ vanishes.

The set of abelian relations  on $U$ (resp.  the set of   germs of abelian relation at a point   $m $ of $M_0$)    has the structure of a  finite dimensional   complex vector space, whose dimension is called the  rank  of the web on $U$ (resp. at $m$).   If the web is   in strong  general position, H\énaut proved in [He2] that   its rank    at a point   does not depend on this point:   abelian relations have then the structure of a local system of coefficients\footnote{In  [He3], he  generalized   this result   in higher codimension.}. When    we    require only the web to be   {   in weak general position}, we shall call ``rank of the web" the maximum\footnote{Notice that,  in this case, the rank at a point is an upper-semicontinuous function of the point.}      of the rank at each point of $M_0$.   

  When the web is   in strong general position,  its  rank    is always upper-bounded, after Chern ([C]), by the number $\pi(n,d)$ of Castelnuovo  (the  maximum of  the arithmetical  genus of an irreducible non-degenerate algebraic curve of degree $d$ in the $n$-dimensional complex projective space {\bb P}$_n$). On the other hand,      
 the rank of an algebraic $d$-web in {\bb P}$_n$ (i.e. the web whose leaves are the hyperplanes belonging to some algebraic curve $\Gamma$ of degree $d$ in the dual projective space {\bb P}$'_n$) is equal to the  arithmetical genus of $\Gamma$: this is, after duality, a theorem coming back to  Abel (see [CG]).  Therefore, the bound $\pi (n,d)$ is optimal for webs in strong general position.  
 
When $n=2$,  the obstruction for the web to have maximal rank $\pi(2,d)=(d-1)(d-2)/2$ is  the Blaschke curvature, which  has been defined    by Blaschke-Dubourdieu ([B]) for  $d=3$, and generalized independently by Panzani ([P])  and  H\énaut ([He1]) for  any $d\geq 3$.  
 
 The first aim of this paper was  to generalize this curvature to any    $n\geq 2$. But some difficulties  appear which do not exist in dimension 2. In fact, we wish   to find   a holomorphic vector bundle with a holomorphic connection such that the vector space  of   abelian relationsis isomorphic to  the vector space of   holomorphic sections with covariant derivative (so that the curvature of this connection will be  an obstruction for the rank of the web to be maximal). Lemma 2-2   below   implies that the  projection 
$R_{k_0-2}\to R_{k_0-3}$ from the space of formal abelian relations at order $k_0-1$ into the space of formal abelian relations at order $k_0-2$ be an isomorphism of holomorphic vector bundles, at least above some everywhere dense open subset $M_0\setminus S$ of $M_0$. Therefore :

 - on one hand ,  $d$ must be equal,    for some $k_0\geq 2$, to the dimension $c(n,k_0)$    of the  vector space of all homogeneous polynomials of degree $k_0$ in $n$ variables with scalar coefficients), 

 - on the other  hand,    some analytical subset $S$ of $M_0$ attached to the web  must have a dimension at most $n-1$ or be empty (this assumption,  is in fact  generically satisfied ; the web is then said to be {\bf ordinary}).   If, exceptionally, $S$ has dimension $n$, the web is said to be {\bf  extraordinary}.

We proved,  by the way, that   whatever be $d$ $\bigl($not necessarily equal to some  $c(n,k_0  )\bigr)$, the rank of  all ordinary  $d$-webs   
 is at most equal to some bound $\pi'(n,d)$ which, for $n\geq 3$, is strictly smaller than the bound $\pi(n,d)$ of Castelnuovo.

   For any $d,\ (d>n)$, denote by $k_0$ the integer ($\geq 1$) such that $c(n,k_0 )\leq d<c(n,k_0+1)$, and set : 
  $$\begin{matrix} \pi'(n,d) &=&0 &\hbox{ when } d< c(n,2) , \ (k_0=1), \\
 &&&\\
 
 &=& \sum_{h=1}^{k_0   }\bigl(d-c(n,h)\bigr) & \hbox{ when }  d\geq c(n,2), \     (k_0\geq  2)  .\\
 
 \end{matrix}$$
 The main results of this paper  are then the two following :
 \begin {theorem} The rank of any ordinary $d$-web    on some $n$-dimensional manifold
 $M_0$ is at most equal to $\pi'(n,d)$. This bound is optimal.  
 \end{theorem}

  \begin {theorem} If $d=c(n,k_0)$, and if the $d$-web is ordinary, the space ${\cal E}= R_{k_0-3} $   of formal abelian relations at order $k_0-2$  is  a holomorphic  vector bundle    of rank $\pi'(n,d)$ over $M_0\setminus S$. There is    a holomorphic connection $\na $ on ${\cal E}$, such that the map $u\mapsto j^{k_0-2}u$ defines an isomorphism from the   space of germs of abelian relations at a point of  $M_0\setminus S$ onto  the vector space of germs of sections of ${\cal E}$ with vanishing  covariant derivative:  $\na ( j^{k_0-2}u)=0$. The curvature of this connection is therefore an obstruction for the rank of the web to reach the value $\pi'(n,d)$.   
 \end{theorem}

For $n=2$,  it happens  \hb 
- that $ S$ is always empty, so that all webs are ordinary, 
\hb - that the  upper-bounds $ \pi(2,d)$ and $\pi '(2,d)$   coincide, 
\hb - that any $d,d\geq 3, $ may be written $d= c(2,k_0 )$, with $k_0=d-1$.\hb 
Thus, we recover the results given locally in [P][He1].

 For $n\geq 3$, the things are not so simple : $S$ may be non-empty,    $\pi'(n,d)$ is strictly smaller than
   $\pi(n,d)$, and not any $d$ may be written $c(n,h)$. Moreover,
  the fact that $dim\ S$
  must be equal to $n$ for the rank of the web to reach the bound $\pi(n,d)$ is   a criterium of practical interest.
  
   In the last section, we give a second  proof of theorem 1-1  in the particular case of webs which are  in  strong general position. We prove also     that all    ordinary affine $d$-webs  in dimension $n$  have rank $\pi'(n,d)$ (hence the optimality of this bound). We give finally significative  examples of the various possible situations.

 We thank A. H\énaut for very helpful conversations.

 \section{Notations and backgrounds }  
 
 Most of the results in this section are well known or easy to prove, so that we shall omit their proof. 
 \subsection{Some algebraic notations}
    
    Let  {\bb R}$_h[X_1,\cdots,X_n]$ denote the vector space of homogeneous polynomials of degree $h$ in $n$ variables, with scalar coefficients (in fact, the field of scalars does not matter). Denote by $c(n,h)$ its dimension $\begin{pmatrix}
	h+n-1\\
	h\\
\end{pmatrix}$. 
   The monomials   $X^L=(X_1)^{\lambda_1}.\cdots(X_n)^{\lambda_n}$ make a basis,   indexed by the set ${\cal P}(n,h)$ of the partitions $L=(\lambda_1,\lambda_2,\cdots,
   \lambda_n)$ of $h$,  with $0\leq \lambda_i\leq h$ for any $i=1,\cdots , n$, and $\sum_{i=1}^n\lambda_i=h$. The number $h$ will be also denoted by $|L|$, and is called the {\emph height} of $L$.
   
  Denoting by $(a)^+$ the number $sup\ (a,0)$ for any real number   $a$,  remember  ([GH]) that the number $$\pi(n,d)=\sum_{h\geq 1}\bigl(d-h(n-1)-1\bigr)^+,$$
   called \emph{the bound of Castelnuovo}, is the maximum of the arithmetical genus of an algebraic curve of degree $d$ in the complex projective space {\bb P}$_n$. 
   
    Define 
    $$\begin{matrix} \pi'(n,d) &=&0 &\hbox{ when } d< c(n,2) , \\
 &&&\\
 
 &=& \sum_{h\geq 1} \bigl(d-c(n,h)\bigr)^+& \hbox{ when }  d\geq c(n,2)   .\\
 
 \end{matrix}$$
   
   \begin{lemma}\hb 
   
   $(i]$ For $n\geq 3$, the inequality $\pi'(n,d)<\pi(n,d)$ holds.
   
   $(ii)$ The  equality $\pi'(2,d)=\pi(2,d)$ holds.
   
     \end {lemma}

  \subsection{Connections adapted to a differential operator :}
  
  Let $E\to V$ be a holomorphic vector bundle on a holomorphic manifold $V$. 
  Remember  the exact sequence of holomorphic vector bundles 
  $$0\to T^*V\otimes E\to J^1E\to E\to 0 ,$$
  where $J^1E$ denotes the holomorphic bundle of 1-jets of holomorphic sections of $E$, and $T^*V$ the holomorphic bundle 
  of 1-forms of type $(1,0)$ on $V$. 
  A $C^\infty$-connection  $\na$ on $E$ of type $(1,0)$ (resp. a holomorphic connection if any)    is a $C^\infty$-splitting (resp. a holomorphic splitting if any) of this exact sequence :
  $$0\to T^*V\otimes E\buildrel{\buildrel \beta\over \longleftarrow }\over \longrightarrow J^1E\buildrel{\buildrel \alpha \over \longleftarrow }\over \longrightarrow E \to  0  .$$
    If $u$ is a $C^\infty$ or holomorphic section of $E$,    its covariant derivative  is given by 
    $$ \na u =\beta(j^1u) .$$
 A  holomorphic linear  differential operator of order $p$ on the holomorphic manifold $V$ is a morphism
 of holomorphic vector bundles $D: J^pE\to F$, for some vector bundles $E$ and $F$ on $V$, to which we associate 
 the map ${\cal D}:u\mapsto D(j^pu)$ from the sections of $E$ into the sections of $F$.  The kernel $R$ of $D$ is a subset of $J^pE$, which is the set of the formal solutions at order $p$ of the equation ${\cal D}u=0$. A  sections $u$ of $E$ such that 
${\cal D}u=0$ is called a "solution" of this equation 

The morphism $D$ induces a morphism $j^1D:J^1(J^pE)\to J^1F$ whose restriction $\widetilde D:J^{p+1}E\to J^1F$ to the sub-bundle $J^{p+1}E\subset J^1(J^pE)$  is called the 
 first prolongation  of $D$,  and the sections $u$ of $E$ such that $\widetilde D(j^{p+1}u)=0$ are also the solutions of the equation 
${\cal D}u=0$ . 

Assume that $R$ is a holomorphic vector bundle above $V$. Then $J^1R$ is a sub-vector bundle of $J^1(J^pE)$, as well as $J^{p+1}E$. The kernel $R'$ of the morphism $\widetilde D:J^{p+1}E\to J^1F$ (the space of formal solutions  at order $p+1$ of the equation ${\cal D}u=0$) is then equal to the 		intersection $J^1R\cap J^{p+1}E$ in  $J^1(J^pE)$. 

 $$\begin{matrix}
 J^1(J^pE)&  & =& & J^1(J^pE)&\buildrel {j^1 D}\over\longrightarrow  &J^1F    \\ 
&&&&&&\\
\cup \uparrow \  \  \  && &&\cup \uparrow \  \  \  & & =  \updownarrow \  \  \    \\
&&&&&&\\
J^1R& \hookleftarrow&R' &\hookrightarrow&J^{p+1}E&\buildrel {\widetilde D}\over\longrightarrow & J^1F    \\
\downarrow & &\downarrow &&\downarrow & & \downarrow     \\
R & =&R &\hookrightarrow&J^pE&\buildrel D\over\longrightarrow & F    \\
 & &\downarrow &&& &     \\
& &V &&& &     \\
   \end{matrix}$$

 \begin{lemma}\footnote{This  lemma, following from the Spencer-Goldschmidt theory ([S]), is used   by H\énaut in [He1]  under an equivalent form.}  The two following assertions are equivalent : 

$(i)$ The projection $R'\to R$ is an isomorphism of holomorphic vector bundles.

$(ii)$ There is a holomorphic connection $\na$ on $R\to V$ such that the map $u\mapsto j^pu$ is an isomorphism from the space  of solutions $u$ of the equation ${\cal D}u=0$ into  the space of sections of $R$ with vanishing covariant derivative. 

Moreover, such a connection $\na $ is unique and defined as being the map $\alpha:R\to J^1R$ equal to the composition of the inclusion $R'\ \hookrightarrow  J^1R$ with the inverse isomorphism $R\to R'$. 
 \end{lemma}
   \n The  connection $\na $ above on $R$   is  said  to be ``completely adapted" to the differential operator ${\cal D}$.

 \subsection{The differential operator for abelian relations }  
  
  Recall that a $d$-web on a  
  $n$-dimensional holomorphic manifold $M$ (in weak general position, with $d>n$) is defined by a $n$-dimensional analytical subspace $W$ of the projectivized cotangent space {\bb P}$(T^*M)$,  on the smooth part of which  the canonical contact form induces a foliation $\widetilde{\cal F}]$. Let $W_0$ be the set of points in $W$ over $M_0$ and $\pi_W:W_0\to M_0$ be the  $d$-fold corresponding covering.  For any complex holomorphic vector bundle $E$ of rank $r$ over $W_0$, let $\pi_*E$ be the holomorphic bundle of rank $d\times r$ over  $M_0$, whose fiber at a point $m\in M_0$ is defined by 
  $$(\pi_*E)_m=\bigoplus _{\widetilde m\in (\pi_W)^{-1}(m)}E_{\widetilde m}\ .$$ 
  The sheaf of holomorphic sections of $\pi_*E$ is then the direct image   of the sheaf of holomorphic sections of $E$ by $\pi_W$. For instance,  an element $\omega$ of $\pi_*(T^*\widetilde{\cal F})$ at $m$  is a family $( \omega_i)$   of forms $\omega_i$ at $m$ ($1\leq i\leq d$),  such that the  kernel of $\omega_i$ contains the
   tangent space $T_i$ to the local foliation ${\cal F}_i$ (and  is therefore equal to it when $\omega_i$
    is not zero).
  
   Let $Tr:\pi_*(T^*\widetilde{\cal F})\to T^*M_0$ be the morphism of holomorphic vector bundles   given    by 
   $$Tr \ \omega \ =\ \sum_{i=1}^d\omega_i.$$ Since the web is in weak general position, it is easy to check that $Tr $ has constant rank $n$.   Then, we
define  a  holomorphic vector bundle   $A$
  of rank    $d-n$ over $M_0$   as the kernel of this morphism.

 \begin{definition} The vector bundle $A$ of rank $d-n$ so defined will be 
   called  \emph{ the Blaschke bundle  } of the web.
   \end{definition} 
   
  Let us  define  similarly    a holomorphic vector bundle $B$ over $ M_0$, of
rank $(d-1) n(n-1)  /2$ as the kernel of the morphism $Tr:\pi_*(\bigwedge^2T^*W_0)\to \bigwedge^2T^*M_0$   given    by 
   $Tr \ \varpi \ =\ \sum_{i=1}^d\varpi_i$. 
    
   \begin{definition}   We 
   call \emph{abelian relation} of the web any
    holomorphic section    $u$ of the Blaschke bundle $A$,  which is  solution 
   of the equation ${\cal D}u=0$, the map ${\cal D}$ denoting the linear first order  
   differential operator $ (\omega_i)\mapsto ( d \omega_i)$ from $A$ to $B$.
   \end{definition}

    The differential operator ${\cal D}$ may still be seen as a linear
morphism
   $D:J^1A\to B$, and the kernel $R_0$ of this morphism is the space of ``formal
   abelian relations at order one".   Hence, a necessary condition for an abelian relation to exist  above an open subset $U$ of
   $M_0$ is that $U$ belongs to the image of $R_0$ by the projection $J^1A\to A$.
   We shall see that it generically not the case for $d<c(n,2)$.

   More generally, denote by  $D_k:J^{k+1}A\to J^kB$    the
   $k^{th}$ prolongation of $D$. The kernel $R_k$ of the morphism $D_k$ is the
   space of  formal
   abelian relations at order $k+1$. Abelian relations may still be seen as the
   holomorphic
   sections $u$ of $A$   such that $j^{k+1}u$ be in $R_k$. Let $\pi_{k+1}$   denote the natural projection $R_{k+1}\to R_k$. 
   Let $\sigma_{k+1 }:S^{k+2}(T^*M_0)\otimes A\to S^{k+1 }(T^*M_0)\otimes B$ be the symbol of $D_{k+1}$, $g_{k+1}$ be the kernel of this symbol, and $K_{k}$  its cokernel. After  snake's lemma, there is a natural map $\partial_k:R_k\to K_k$ in such a way that we get  a commutative diagram\footnote{Notice that the $R_k$'s, $g_k$'s and $K_k$'s are not necessarily vector bundles : exactitude has   to be understood as the exactitude in  the corresponding diagrams on the fibers at any point of $M_0$. }, with all lines and columns exact:
  $$\begin{matrix}
 & & 0 & &0 & & 0 & &R_k &     \\ 
 & & \downarrow  & & \downarrow  & & \downarrow   & &  \ \ \ \   \downarrow   \partial_k &     \\
  
 0 & \rightarrow & { {g}}_{k+1}& \rightarrow & S^{k+2}(T^*M_0)\otimes A &\buildrel{\sigma_{k+1}}\over\longrightarrow & S^{k+1}(T^*M_0)\otimes B & \rightarrow &  K_k &  \rightarrow  0    \\
 
 & & \downarrow  & & \downarrow  & & \downarrow   & & \downarrow  &     \\
  
 0&\rightarrow &R_{k+1} &\rightarrow &J^{k+2}A &\buildrel{D_{k+1}}\over\longrightarrow  &J^{k+1}B & \rightarrow & coker \ D_{k+1}  &  \rightarrow  0      \\
 
 & &\ \ \ \  \downarrow\pi_{k+1}  & & \downarrow  & & \downarrow   & & \downarrow &     \\
 
  0&\rightarrow &R_{k } &\rightarrow &J^{k+1}A &\buildrel{D_{k }}\over \longrightarrow &J^{k }B & \rightarrow & coker \ D_{k }&    \rightarrow  0       \\
  
  & &\ \ \ \  \downarrow  \partial_k & & \downarrow  & & \downarrow   & &  \downarrow &     \\
  
  & &K_k& &0 & & 0 & &0 &     \\

  \end{matrix}$$
   In this diagram, we  allow $k$ to take the value $-1$, with the convention $R_{-1}= A$, $J^{-1}B=0$. 
   \n In the sequel,  
  
   the index
  $i $ will run from $1$ to $d$, 
  
  the indices $\alpha,\beta,\cdots$ will run from $1$ to $n-1$, 
  
  and 
  the indices $ \lambda,\mu,\cdots$ will run from $1$ to $n$. 
  
  For any  holomorphic function $a$ and local coordinates $x$, $a '_\lambda$ will denote  the partial derivative  $a '_\lambda = {{\partial a }\over{\partial x_\lambda}}$. More generally, 
  for any  partition  $L=(\lambda_1,\lambda_2,\cdots,\lambda_n)$ of    $ |L|=\sum_\mu\lambda_\mu$,  ($ L\in {\cal P}(n,|L|)$ ), we denote by   $(a)'_{_L} $
      the corresponding higher order derivative 
   $ (a)'_{_L} = {{\partial^{|L|} a }\over{(\partial x_1)^{\lambda_1} \cdots(\partial  x_n)^{\lambda_n}}}$.
   
   \section {Definition   of the analytical set $S$}
   This definition is different according to the inequalities $d<c(n,2)$ or $d\geq c(n,2)$. However, in both cases, this
     set will satisfy to the  \begin{lemma} The set $S$ is an analytical set which has generically a dimension $\leq n-1$ or is empty.
 \end{lemma}

     \subsection{Definition of $S$ in the case $n<d<c(n,2)$:}

     For $n<d<c(n,2)$,   $S$ will denote  the subset of elements $m\in M_0$ such that the vector space\break  $(R_0)_{_m}=(\pi_{0})^{-1}(A_m)$ has dimension at least 1. 
     The proof of lemma 3-1 for $d<c(n,2)$ will be given in the next section, after the description of the map $\pi_{0}$.
\n {\bf Proof of theorem 1-1 in the case $d<c(n,2)$ :}   {\it 
For  $d<c(n,2)$, all ordinary $d$-webs of codimension one have rank 0.}
 \n  The fact that all ordinary $d$-webs   have rank 0 for  $d<c(n,2)$ on $M_0\setminus S$ is a tautology, because of the definition of $S$. Therefore, they still have rank 0 on all of $M_0$,     because  of the semi-continuity of the rank at a point. 
 
 \rightline{QED}

 \subsection{Definition of $S$ in the case $d\geq c(n,2)$}
 
 Let  $ \eta_i$ be an integrable 1-form defining the local foliation ${\cal F}_i$
 of a $d$-web   of codimension 1. For any integer $h\geq 1$,  let $(\eta _i)^h $ be the $h^{th}$ symetric power of $\eta_i$
 in the space $S^h(T^* M_0)$ of homogeneous polynomials of degree $h$ on the complex tangent tangent bundle  $T M_0$. For any $m\in M_0$, let $r_h(m)$ be the dimension of the subspace $L_h(m) $ generated in $S^h(T^*_mM_0)$ by the $(\eta _i)^h(m)$'s with $1\leq i\leq d$ (not depending on the choice of the $\eta_i$'s).  We have obviously the 
 \begin{lemma}The following inequality holds :
 $$r_h(m)\leq {\rm min} \bigl(d,c(n,h)\bigr).$$
     \end {lemma}
     In particular,  $r_1 \equiv n$, since we assume   $d>n$ and    the web to be { in weak general position}. And $r_h\equiv d$ for $h>k_0$, where $k_0$ denotes the integer such that 
 $$c(n,k_0)\leq d<c(n,k_0+1).$$ 
\n {\bf Definition of $S$ for $d\geq c(n,2)$  (i.e. $k_0\geq 2$) :} Let $S_h$ be the set of points $m\in M_0$ such that $r_h(m)<{\rm min} \bigl(d,c(n,h)\bigr)$, and set : $S=\bigcup_{h= 2}^{k_0}S_h$.

 \n {\bf Proof of lemma 3-1 for $d\geq c(n,2)$ : } Let   $(x_\lambda)_{1\leq \lambda\leq n}$ be a  system of holomorphic  local coordinates near  a point   $m_0$  of $M_0$, and assume that the local foliation ${\cal F}_i$ is defined by \hb $ \sum_\lambda p_{i\lambda}dx_\lambda=0$. Then 
 $r_h(m)$ is the rank of the matrix $$P_h=(\!(C_{i\ L}^{(h)})\!)_{_{i, L}}$$ of size $d\times c(n,h)$ at point $m$,  where   $ 1\leq i\leq d$,  and  $L$ runs through the set ${\cal P}(n,h)$  of the partitions \break $L=(\lambda_1,\lambda_2,\cdots ,\lambda_n)$ of $h$ (i.e.  $\sum_{s=1}^n\lambda_s=h $), and where  
   $C_{i\ L}^{(h)}=\prod _{s=1}^np_{_{i\ \lambda_s}}  $.  Thus, for $2\leq h\leq k_0$, $S_h$ is locally defined by the vanishing of all determinants of size $c(n,h)$ in $P_h$. Exceptionnaly, it may happen that all  of these determinants are identically 0,  so that $S_h$ will have dimension $n$. But generically,   these determinants will vanish  on hypersurfaces or nowhere. 
   \rightline {QED}

 
  \section {Computation of $R_0$} 
    
   Locally, the $d$-web is defined over $M_0$ by  a family of 
   1-forms $\eta_i=dx_n-\sum_\alpha p_{_{i\alpha}}(x)\ dx_\alpha$, which we still may write 
   $\eta_i= -\sum_\lambda p_{_{i\lambda}}(x)\ dx_\lambda$  
  with the convention $p_{_{in}}\equiv -1$.

   \begin{lemma} The integrability conditions   may   be written locally :
$$ (p_{i\lambda} )'_\mu-(p_{i\mu } )'_\lambda +p_ {i\mu } (p_{i \lambda } )'_n-p_ {i\lambda}(p_{i\mu })'_n \equiv 0 
  \hbox   {\ \ \  for all triples  }  (i, \lambda, \mu)    . $$
 \end{lemma}
  Proof :  Let $\eta = -\sum_\lambda p_{_{ \lambda}}(x)\ dx_\lambda$  be a holomorphic 1-form. Then 
  $$ \eta\wedge d\eta   =    
 \sum_{\lambda <\mu <\nu  }\biggl[ p_\lambda \Bigl((p_\nu )'_\mu -(p_\mu )'_\nu  \Bigr)
 + p_\mu \Bigl((p_\lambda )'_\nu -(p_\nu )'_\lambda  \Bigr)+
  p_\nu \Bigl((p_\mu) '_\lambda -(p_\lambda )'_\mu \Bigr) \biggr] d   x_\lambda  \wedge d   x_\mu  \wedge d x_\nu .$$
 Then, when $p_n\equiv-1$, we  observe   that      the vanishing of all terms in $d x_ \lambda\wedge d x_\mu \wedge d x_n $
   implies the vanishing of all other terms,    hence the lemma.
   
   \rightline{QED}
  
  A section $(\omega_i)_i$ of $A$ is locally given by
   the $d$ functions $f_i$ such that $\omega_i=f_i\ \bigl(\sum_\lambda p_{_{i\lambda}}(x)\
   dx_\lambda\bigr)$  satisfying to the identities
   $$ (E _\lambda ) \hskip 2cm\sum_i\ p_{_{i\lambda}}f_i\equiv 0 \hskip 1cm\hbox{ for any }\lambda , $$
   hence, by derivation,    
   $$ (E _{\lambda,\mu  }) \hskip 1cm\sum_i\ ( p_{_{i\lambda}}f_i)'_\mu \equiv 0 \hskip 1cm\hbox{ for any }\lambda   ,\mu . $$
    
     \begin{lemma} For the family  $(f_i)_i$ to define  an abelian relation, it is necessary and sufficient that the identities be satisfied 
$$( F_{i\alpha} ) \hskip 4cm (f_i)'_ \alpha\equiv  -\bigl(f_i p_{i \alpha} \bigr)'_n  \hbox {\ \  for all pairs }  (i,\alpha).$$ \end{lemma}
  Proof : In fact, a holomorphic 1-form $f\bigl(\sum_\lambda p_{_{ \lambda}}(x)\ dx_\lambda\bigr)$ is closed iff 
  $(fp_\lambda)'_\mu=(fp_\mu )'_\lambda $ for all pairs   $( \lambda, \mu)$   such that    $\lambda <\mu  .$ But, because of the integrability conditions,  it is sufficient that this relation be satisfied when  $\mu=n$, for it to be  satisfied with all other $\mu$'s.
  
  \rightline{QED}
  \begin{lemma} When the family  $(f_i)_i$  defines  an abelian relation, the identities $(E _{\lambda,\mu  })$ and $(E _{ \mu,\lambda  })$ are the same.   \end{lemma}
  Proof : In fact, under the assumption,  $(f_ip_{i\lambda})'_\mu\equiv (f_ip_{i\mu })'_\lambda $ for all pairs   $( \lambda, \mu)$, hence the lemma  by summation with respect to $i$.
  
  \rightline{QED}
  The identity $(f_i)'_ \alpha\equiv  -\bigl(f_i p_{i \alpha} \bigr)'_n $ means that it is sufficient to  know  the $(f_i)'_n$ to know the other partial derivatives $(f_i)'_\alpha $ of a  family  $(f_i)_i$   defining   an abelian relation. 
  
  Hence, writing $w_i=(f_i)'_n$, and combining $(E _{\lambda,\mu  })$ and $( F_{i\alpha} )$, we get :
  \begin{corollary}  The elements of $R_0$ above a given element $(f_i)_i$ in $A$ map bijectively onto the solutions of the linear system $\Sigma_0$. 
  $$(\widetilde E _{\lambda,\mu  }) \hskip 2cm\sum_i\   p_{_{i\lambda}}p_{_{i\mu }}w_i   \equiv 
  \sum_i\   f_i\bigl[(p_{_{i\lambda}})'_\mu-p_{_{i \lambda }}(p_{_{i\lambda}})'_n\bigr]  $$
  of $c(n,2)$ equations $(\widetilde E _{\lambda,\mu  })$ with $d$ unknown $w_i$.
  \end{corollary}
  Notice that the matrix of the system $\Sigma_0$ is the matrix  $P_2=(\!(C_{i\ L}^{(2)})\!)_{_{i, L}}$ seen in the previous section.
 \n  {\bf Proof of the lemma 3-1 in the case $d<(c(n,2)$ :} Generically, the system $\Sigma_0$ has rank $d$.  Hence $S$ is defined locally by the vanishing of all characteristic determinants which are generically not all identically zero. If $\Sigma_0$ has a rank $r$ smaller than $d$, the vanishing of all determinants of size $r+1,\cdots,d$ in the matrix $P_2$ have to be added to the vanishing of all charateristic determinants. 
 
 \rightline{QED}

      \section {Computation of $R_{k }$ $(k\geq 1)$ :}
      
      For  any pair of multi-indices  of derivation $L=(\lambda_1, \cdots,\lambda_s,\cdots,\lambda_n)$, and $H=(h_1,h_2,\cdots,h_n)$, $L+H$ will denote the multi-index 
   $(\lambda_1+h_1,\lambda_2+h_2,\cdots,\lambda_n+h_n)$. We define similarly 
 $L-H $ if $\lambda_\mu\geq h_\mu$ for all $\mu$'s.  For any $\lambda$, $1_\lambda$ will denote the muti-index 
 with all $\lambda_\mu$'s equal to zero for $\mu\neq \lambda$ and $\lambda_\lambda=1$. 
 By definition the height $|L|$ of $L$ is the sum $\sum_s\lambda_s$.

   By  derivation of the identities  $(E _\lambda )$, the elements of $J^kA$ are characterized by the identities  
   $$ (E _{\lambda,L }) \hskip 1cm\sum_i\ ( p_{_{i\lambda}}f_i)'_L\equiv 0 \hbox{ for any $\lambda$ and for any multi-index $L$ of height }|L|\leq k .   $$

   \begin{lemma}
   If $(f_i)_i$ is an abelian relation, the relation $(E _{\lambda,L })$ remains unchanged by permutation of all    the indices of $ L\cup \{\lambda\}$. 
      \end{lemma}
      \n  Proof : The left hand term of this identity is obviously symetric with respect to the indices of $L$. Thus, it is sufficient to prove that the identities $(E _{\lambda, \mu }) $ and $(E _{ \mu,\lambda}) $ are the same, which we know already.   \rightline{QED} 
      
      The identity $(E _{\lambda,L })$ above will now be denoted by $(E _{H})$, where  $H=L+1_\lambda$.

  \begin{lemma}\hb 
  $(i)$ If $(f_i)_i$ is an abelian relation,  all partial derivatives  $ (f_i)'_L $ may be written  as a  linear combination 
  $$(\widetilde F_{i L}) \hskip 3cm\widetilde {(f_i)'_L}\equiv \sum_{k=0}^{|L|}\ D_{i\ L}^{(k)}\ .\ (f_i)' _{n^k}$$ 
  of $f_i$ and of its partial derivatives  
   $ (f_i)' _{n^k}=  {{\partial^{(k)} f_i }\over{(\partial x_n)^{k}  }}$  
  with respect to the only variable $x_n$, with coefficients $D_{i\ L}^{(k)}$ not depending on the $f_i$'s.   
\n 
  $(ii)$ If $L=(\lambda_1, \cdots,\lambda_s,\cdots,\lambda_n)$, the coefficient $D_{i\ L}^{(|L|)}$ of highest order is equal to $(-1)^{|L|}\prod _{s=1}^np_{_{i }}^{\lambda_s} $, i.e. is equal to $(-1)^{|L|}C_{i\ L}^{(\vert L\vert)}$.  
      \end{lemma}
      \n Proof :   We get   the    lemma  by derivation of the identities $( F_{i\alpha} )$ and an obvious induction on the height $|L|$ of $L$. 
      
      \rightline{QED}
      
       The lemma above  means that it is sufficient to know the $ (f_i)' _{n^{k+1}}$ to know the other partial derivatives $(f_i)'_L $ in the $(k+1)$-jet of a  family  $(f_i)_i$   defining   an abelian relation. 
       
  Hence, writing $w_i=(f_i)'_{_{n^{k+1}}}$, and combining $(E _{L }$ and $(\widetilde F_{i L}) $, we get :
  \begin{corollary}  The elements of $R_{k }$ above a given element $ a_0^{(k-1 )} $ in $R_{k-1 }$ map bijectively onto the solutions of a linear system $\Sigma_{k}$ of $c(n,k+2)$ equations $(\widetilde E _{L })$ with $d$ unknown $w_i$  
  $$(\widetilde E _{L }) \hskip 2.7cm\sum_i\   C_{i\ L}^{(k+2)}\ w_i   \equiv \Phi_L\bigl(a_0^{(k-1 )}\bigr)  ,$$
 where $L$ runs through the set ${\cal P}(n,k+2)$ of partitions of $k+2$,   and 
   where the second member $\Phi_L\bigl(a_0^{(k -1)}\bigr)$ depends only on $ a_0^{(k-1 )}\in R_{k-1 }  $.  
In particular, the symbol $\sigma_{k}$ of $D_{k}$ is defined by   the matrix $$P_{k+2}=(\!(C_{i\ L}^{|L|})\!)_{_{i, L}},$$    $ 1\leq i\leq d$, $|L|=k+2 $, of size $d\times c(n,k+2).$ 
   \end{corollary}

  \begin{theorem}  Assume the $d$-web to be ordinary.
    The   map  $\pi_{k } :R_{k  }\to R_{k-1} $ is surjective for $k \leq k_0-2 $ and injective for $k  > k_0-2 $ above the open set $U= M_0\setminus S$ of $M_0$.
     For $k \leq k_0-2 $, $R_{k }$ is a holomorphic bundle of rank 
  $\sum_{h=1}^{k+2  }\bigl(d-c(n,h)\bigr) $  over $M_0\setminus S$.  
 \end{theorem}
   \n Proof : In fact, $P_{k+2}$ is precisely  the matrix of the system $\Sigma_k$. Thus, 
 for $k\leq k_0-2$  and off $S_k$, the space of solutions of $\Sigma_k$ is an affine space of dimension $d-c(n,k+2)$.

 \rightline{QED}
\n {\bf Proof of   theorem 1-1 in the case $d>c(n,2)$ :}
   {\it The rank of a ordinary $d$-web is at most equal to the number 
   $$\pi'(n,d)=\sum_{h=1}^{k_0  }\bigl(d-c(n,h)\bigr).$$}
   
   For $k>k_0-2$, the symbol $\sigma_k$ is necessarily injective off $S$. In fact, the matrix defining this symbol in the corollary above contains the matrix defining $\sigma_{k-1}$ when we choose the coordinates and the forms $\eta_i$'s
   so that $p_{in}\equiv -1$. Hence, since $g_{k_0-2}=0$ off $S$, $g_k=0$   off $S$ for all $k\geq k_0-2$. Consequently, above a given element in $R_{k_0-2}$, there exists at most one infinite jet of abelian relation, hence one germ of abelian relation since the framework is analytic. We deduce that the rank of the web is at most $\pi'(n,d)$ off $S$, hence everywhere (semi-continuity of the rank). We shall see  in the next section that any ordinary affine web in strong  general position has rank $\pi'(n,d)$, hence the optimality. 
   
    \rightline{QED}  
 \n {\bf Proof of   theorem 1-2   :}
  When  $d=c(n,k_0 )$,     ${\cal E}=R_{k_0-3}|_{M_0\setminus S}$ is a vector bundle   of rank    $\pi'(n,d) $, and  $\pi_{k_0-2} : R_{k_0-2  }\to R_{k_0-3} $ is an isomorphism of vector bundles over $M_0\setminus S$. Therefore, we just have to use  lemma 2-2. 
 In particular, for $d=c(n,k_0 )$, the rank of a ordinary $d$-web is  $\pi' (n,d) $ iff the curvature of the previous connexion vanishes.
 \rightline{QED}

    \section{ Examples }
     \subsection{ Case $n=2$ :}
    
      We recover the results of [P][He1]).  In fact, in this case  :
   
- the set $S$ is   always empty (all determinants occuring in the computation of the symbols $\sigma_k$ are determinants of Van-der-Monde for $k\leq k_0-2$,   vanishing nowhere on $M_0$); thus, all webs are ordinary. 
   
 -  any $ d $ is equal to $c(2,d-1)$, 
   
 -   the rank  $\sum_{h=1}^{d-2}\bigl(c(2,d-1)-c(2,h)\bigr)$ of $R_{d-4}$ is   equal to $(d-1)(d-2)/2$.

\subsection{ Case  $n=3$, $d=6$ :}

Use    coordinates $ x,y,z $ on  {\bb C}$^3$, with $n=3$, and  $c(n,2)=6$.  
Let $a,b, c,e, h$ be five  distinct  complex numbers, all  different of $0 $, and let $ \psi $ be some holomorphic function   of $y$. Let $W$ be the 6-web of codimension  1 on 
{\bb C}$^3$ defined by the 1-forms   $\eta_i=dz-p_idx-q_idy,\ 1\leq i\leq 6$, where :\break   
$(p_1,q_1)= (0,0)\ ,\  (p_2,q_2)= (a,a^2)\ ,\   (p_3,q_3)= (b,b^2)$,   $(p_4,q_4)= (c,c^2)\ ,\  (p_5,q_5)=  (e,e^2)$   and   $ 
(p_6,q_6)= (h, \psi)$.  
The   system $\Sigma_0$ is then equivalent to 
    $$\begin{pmatrix}
 a&b&c&e& \\ 
 a^2&b^2&c^2&e^2 \\
 a^3&b^3&c^3&e^3& \\
 a^4&b^4&c^4&e^4& \\ 
\end{pmatrix}\begin{pmatrix} w_2\\w_3\\w_4\\w_5   \end{pmatrix}=
\begin{pmatrix}   0\\0\\0\\f_6.\psi'\end{pmatrix},  $$
 to which we add the equations $(h^2-\psi)w_6=0$, and $w_1+w_2+\cdots +w_6=0$. 
The system is a system of Cramer off the  locus $S$ which is 

- the surface    of equation $ \psi(  y ) =h^2$ in general,  

- the empty set when $\psi $ is any constant different of $ h^2$, 

- all of  {\bb C}$^3$ (the extraordinary  case)  for $ \psi\equiv   h^2$.

 Let's precise the   connection and its curvature for $ \psi\neq   h^2$ in the two first cases (the ordinary case).    

A section $(f_i)$ of $A$ is defined by $(f_4,f_5,f_6)$, since $(f_1,f_2,f_3)$ can be deduced using 
 the equations $\sum_if_i=0$,  $\sum_ip_if_i=0$   and $\sum_iq_if_i=0$. 
Setting $\Delta= abce(b-a)(c-a)(e-a)(c-b)(e-b)(e-c)(h^2-\psi)$,  $K={{\psi'}\over{\Delta(h^2-\psi)}}$, $\Delta_4=abe((b-a)(e-a)(e-b)$,  and  $\Delta_5=abc((b-a)(c-a)(c-b)$,we get :\hb  $ w_4=-K\Delta_4,w_5= K\Delta_5,w_6=0$, hence $ u_4=cK\Delta_4,  u_5=-eK\Delta_5,u_6=0$ and 
$ v_4=c^2K\Delta_4,$\break $u_5=-e^2K\Delta_5,\ u_6=0$. With respect to the trivialization $(\sigma_4,\sigma_5,\sigma_6)$ of $A$ given by \hb 
$\sigma_4=(f_4\equiv 1,f_5\equiv 0, f_6\equiv 0)$, $ \sigma_5=(f_4\equiv 0,f_5\equiv 1, f_6\equiv 0)$ and 
$\sigma_6=(f_4\equiv 0,f_5\equiv 0, f_6\equiv 1)$, the matrix of the connection is : 
 $$\omega={{\psi'}\over{\Delta(h^2-\psi)}}\begin{pmatrix} 0&0&\Delta_4\ \eta_4\\ 0&0&-\Delta_5\  \eta_5\\0&0&0\\    \end{pmatrix} ,  $$ 
hence the  curvature  $$ 
  \Omega= {1\over \Delta} \Bigl({{\psi'}\over{  h^2-\psi }}\Bigr)'\begin{pmatrix} 0&0&\Delta_4(dy\wedge dz+c\ dx\wedge dy)\\ 0&0&- \Delta_5(dy\wedge dz+e\ dx\wedge dy)\\0&0&0\\    \end{pmatrix}  .$$ 
  We observe that $\sigma_4$ and $\sigma_5$ are linearly independant abelian relations. We knew it already since the first integrals 
  $(z-cx-c^2y)$ and $(z-ex-e^2y)$ of ${\cal F}_4$ and ${\cal F}_5$ respectively are linear combinations of 
  the first integrals $z$, 
  $(z-ax-a^2y)$ and $(z-bx-b^2y)$  of ${\cal F}_1$, ${\cal F}_2$ and ${\cal F}_3$. Thus, when $ h^2-\psi$ does not vanish,  the rank of the 6-web is at least 2, and has the maximum possible value $\pi'(3,6)=3$ in the ordinary case   if and only if  ${{\psi'}\over{  h^2-\psi }}$ is constant, 
  that is if there exists two scalar constant $C$ and $D$ $(C\neq 0)$, such that $$\psi(y)=h^2+Ce^{Dy},$$ in particular 
   for $\psi=$ constant (case $D=0$). For given $C$ and $D$ as above, $K=-D/\Delta$, and $\sigma_6-K\Delta_4\sigma_4
 +K\Delta_5\sigma_5$ is an abelian relation (the function $\psi$ occurs in the computation of $f_1$, $f_2$ and $f_3$ from $f_4=-K\Delta_4$, $f_5= K\Delta_5$ and $f_6=1$).  
 
  The extraordinary case $\psi\equiv h^2$ will be seen in the next subsection.

  \subsection{Ordinary affine webs and optimality of the bound $\pi'(n,d)$ :}
 For any pair $(n,d)$ with $d>n$,  give $d$ linear forms $(l_1,l_2,\cdots,l_d)$.   Let ${\cal F}_i$ be  the foliation defined in {\bb C}$^n$ by the parallel hyperplanes $l_i=constant $, and $W$ be  the $d$-web defined by these $d$ foliations.   Let $k_0$ be the integer such that $c(n,k_0)\leq d<c(n,k_0+1)$. Let {\bb P}'$_{n-1}$ denote the hyperplane at infinity
      of    the $n$-dimensional projective space {\bb P}$_n =\ ${\bb C}$^n \coprod $\ {\bb P}$_{n-1}$ :  the parallel hyperplanes of the pencil $l_i=constant  $
      meet at infinity along a hyperplane  of  {\bb P}$_{n-1}$, i.e. define an  element $[l_i]$ of the dual projective space {\bb P}$'_{n-1}$\ .
      \begin{remark} {\rm  The web $W$ on {\bb C}$^n$ extends to a web on {\bb P}$_n$ (with singularities) which is 
      algebraic.  It is in fact dual to the union of $d$ straight lines in the dual projective space {\bb P}$'_{n}$ .   }
            \end{remark} 
            \begin{definition}
            Such a  web on {\bb P}$_{n}$ will be said an {\bf affine $d$-web}, and it is said ordinary if its restriction to 
      {\bb C}$^{n}$ is ordinary.\end{definition}
      \begin{lemma}
      The two following properties are equivalent :
      
     $ (i)$ The affine $d$-web above is ordinary.
      
    $(ii) $ For any $h$,  $(1\leq h\leq k_0)$, there exists $c(n,h)$ points   among the $d$ points $[l_i]$,  which do not belong to a same algebraic $($reducible or not$)$ hypersurface of degree $h$   in {\bb P}$'_{n-1}$. 
      
      \end{lemma}
  Proof :  For $d\geq c(n,h)$, $(h\geq 1)$, assume that  the matrix $(\!(C_{i\ L}^{(h)})\!)_{_{i, L}}$,    $ 1\leq i\leq d$, $|L|=h  $, of size $d\times c(n,h )$ has rank $<c(n,h)$. This means that the determinant of any  square sub-matrix of size $c(n,h)$   vanishes. Saying that the determinant of  the sub-matrix given for instance  by the $c(n,h)$ first $l_i$'s vanishes means precisely that  $[l_1],[l_2],\cdots,[l_{c(n,h)}]$ belong to some hypersurface of degree $h$ in {\bb P}$'_{n-1}$ (may be reducible), and same thing for any other subset of $c(n,h)$ indices $i$. Hence the lemma.  
     
     \rightline {QED}
          
      \begin{theorem}\hb\hb 
   $(i)$ If    the $d$ points $[l_i]$ are in general position in {\bb P}$'_{n-1}$ $($i.e. if any $n$ of the $d$ linear forms  $ l_i $ are linearly independant$)$,    the above affine $d$-web        has a rank $\geq \pi'(n,d)$. 
     
    \n $(ii)$  It has  exactly    rank
     $\pi'(n,d)$  iff it  is ordinary.

            \end{theorem}
            Proof : For any $h$, consider the vector space $L_h$ generated by the $h^{th}$ symetric products $(l_i)^h$ of the $l_i$'s, and denote by $r_h$ the dimension of this   vector space.  
 The rank of an affine  web in strong general position is $\sum_{h=1}^{k_0}(d-r_h)  $ (see Tr\épreau   [T], section 2).  Hence we have only to prove that  $r_h\leq c(n,h)$ in general, and $r_h= c(n,h)$ in the ordinary case.  But  this is obvious   after the previous lemma, since the dimension of $L_h$ is exactly the rank of the matrix   $P_h=(\!(C_{iL}^{(h)})\!)_{_{i, L}}$ with $i\leq d, |L|=h $ in the linear system $\Sigma_{h-2}$ of corollary 5-3.

\rightline{QED} 
   \begin{corollary} The bound $\pi'(n,d)$ for the rank of a ordinary web is optimal. 
            \end{corollary} 
  Proof : Remember that a algebraic hypersurface of degree $h$ in {\bb P}$'_{n-1}$ is defined in general by the data of $c(n,h)-1$ of its points : thus, the property (ii) of lemma  6-3 above is generically satisfied, so that    there exists ordinary affine $d$-webs in dimension $n$ for any $(n,d)$. 

\rightline{QED}   
For instance,  an affine 6-web  on {\bb P}$_3$  in  strong general position   will be extraordinary of rank 
$4\ \bigl(=\pi(3,6)\bigr)$    or ordinary of rank $3\ \bigl(=\pi'(3,6)\bigr)$,  according to the fact that the six points $[l_i]$ belong or not to a same conic\footnote{In the extraordinary case, it is amazing to deduce    the fourth abelian relation from the theorem of the hexagon (Pascal).} of {\bb P}$'_2$ : when $\psi$ is constant in the example $(p_1,q_1)= (0,0)\ ,\  (p_2,q_2)= (a,a^2)\ ,\   (p_3,q_3)= (b,b^2)$,  $(p_4,q_4)= (c,c^2)\ ,\  (p_5,q_5)=  (e,e^2)$   and   $ 
(p_6,q_6)= (h,   \psi)$ of the previous subsection, the five first points belong to the conic 
 of  equation $q=p^2$, hence the dichotomy according to the fact that   $\psi$ is equal or different of $h^2$. 
   
\n {\bf  New  proof of the theorem 1-1 for $d\geq c(n,2)$  in the particular  case of ordinary webs which are  in  strong general position :} 
  
  The rank of a web       in strong general position  is upper-bounded by the rank of the ``tangent affine web" at a point (see  again  [T]). Since the  ordinaryity of a    web near a point is equivalent to the ordinaryity of the tangent affine web at that point,  theorem 1-1 is also  a corollary of the previous theorem 6-4,   when  the ordinary      web   is   in strong general position.  
  
 \rightline{QED}

\subsection{Another  example of ordinary web of maximal rank 26 for  $n=3$, $d=15$ }     
The following  example has been given to us by J.V.  Pereira and L. Pirio  ([PP]).  

In {\bb    C}$^3$ with coordinates $x,y,z$, take the 15-web defined by  the  ten pencils of planes, 4 pencils of quadratic cones or cylinders, and a last  pencil of quadrics,  respectively defined by  the first integals $u_i$ ($1\leq i\leq 15$) :

\n $u_1=x$, $u_2=y$, $u_3 =z$, 
 $u_4=\frac{z}{z-y}$,  $u_5=\frac{z}{z-x}$,  $u_6=\frac{y}{y-x}$, 
 $u_7=\frac{1-z}{y-z}$, 
$u_8=\frac{1-z}{z-x}$, $u_9=\frac{1-y}{y-x}$, 
  $u_{10}=\frac{x-y}{z-y}$, 

\n $u_{11}=\frac{z(1-y)}{z-y}$, $u_{12}=\frac{z(1-x)}{z-x}$, $u_{13}=\frac{y(1-x)}{y-x}$, $u_{14}=\frac{z(x-y)}{x(z-y)}$, 

\n $u_{15}=\frac{(1-z)(y-x)}{(1-x)(y-z)}$.

Denote by $m_X=[1;0;0;0]$, $m_Y=[0;1;0;0]$, $m_Z=[0;0;1;0]$ and  $m_T=[0;0;0;1]$ the edges of the standard projective frame, and by $\Omega_X=[0;1;1;1]$, $\Omega_Y=[1;0;1;1]$, $\Omega_Z=[1;1;0;1]$ and $\Omega_T=[1;1;1;0]$ the barycenters of the faces of the previous tetrahedron.  The foliations 7, 8, 9 and 10 are respectively the foliations given by the pencil of planes through the line $m_X\Omega_X$ (resp. $m_Y\Omega_Y$, $m_Z\Omega_Z$ and $m_T\Omega_T$. 

 The foliation  11 is  the pencil  of the quadratic cylinders of summit 
$m_X$ having for basis the conics through $\Omega_X,m_Y,m_Z,m_T$ in the plane $X=0$. We get similarly the foliations 12, 13, and 14 by permutation of the letters $X,Y,Z, T$. 

The 5-subweb $(1,2,3,10,14)$ is then the 5-web of cones of summit $m_T$ over the  Bol's 5-web defined by the 4 points $m_X,m_Y,m_Z, \Omega_T$ in the plane  at infinity  $T=0$.  
We can do the same with the other 5-subwebs $(2,3,4,7,11)$, $(1,3,5,8,12)$ and $(1,2,6,9,13)$.

L. Pirio checked, using Mapple,  that this 15-web is ordinary. Since $15=c(3,4)$, this web has a curvature, and we   could theoretically check that this curvature vanishes. In fact,  Pereira and Pirio exhibited directly 26 independant abelian relations : since $\pi'(3,15)=26$, the rank of the web is exactly 26 (while $\pi(3,15)=42$).

  \n {\bf References}

 \n [B] W. Blaschke, Einf\"uhrung in die Geometrie der Webe, Birkha\"user, Basel, 1955. 
 \n [C] S.S.Chern, Abz\"ahlungen f\"ur Gewebe, Abh. Hamburg 11, 1936, 163-170.  
\n [CG] S.S.Chern and P.A. Griffiths, Abels theorems and webs, Jahr. Deutsch Math. Ver.  80, 1978, 13-110, and 83, 1981, 78-83. 
\n [He1] A. H\énaut, On planar web geometry through abelian relations and connections, Annals of Mathematics 159, 2004, 425-445.
\n [He2] A. H\énaut, Syst\èmes diff\érentiels, nombre de Castelnuovo, et rang des tissus de {\bb C}$^n$, Publ. RIMS, Kyoto University, 31(4), 1995, 703-720.
\n [He3] A. H\énaut, Formes diff\érentielles ab\éliennes, bornes de Castelnuovo et g\éom\étrie des tissus, Comment. Math. Helv. 79, 2004, 25-57.
\n [GH] P.A.  Griffiths and J. Harris, Principles of algebraic geometry, John Wiley \& Sons, New York, 1978. 
\n [P] A. Pantazi, Sur la d\étermination du rang d'un tissu plan, C.R. Acad. Sc. Roumanie 2, 1938, 108-111.
\n [PP] J.V. Pereira and L. Pirio : private conversation. 
\n [S] D.C. Spencer, Selecta, 3, World Sci. Publishing  Co. Philadelphia, 1985.
\n[T] J.M. Tr\épreau, Alg\ébrisation des tissus de codimension 1, la g\'en\'eralisation d'un   th\éor\ème de Bol, Inspired by S.S. Chern,     edited by  P. Griffiths, Nankai Tracts in Mathematics, vol. 11, World Scientific and Imperial College Press, 2006. 

 \bigskip
\n Vincent Cavalier, Daniel Lehmann\hb D\'epartement des Sciences Math\'ematiques, CP 051,
Universit\'e de Montpellier II \hb Place Eug\`ene Bataillon, F-34095 Montpellier
Cedex 5, France\hb   cavalier@math.univ-montp2.fr;\ \ \ lehmann@math.univ-montp2.fr

  \end{document}